\documentclass[11.5pt]{article}

  \usepackage{amsmath}
   \allowdisplaybreaks  

\def\R{{\mathbb{R}}}

\def\loc{{\mathrm{loc}}}

\usepackage{graphicx,amssymb,mathrsfs, latexsym,amsfonts,titlesec,amscd,epsfig,cite
                                                                                   }
\usepackage{indentfirst} 
 \usepackage{enumerate}
\usepackage{units} 

\DeclareMathOperator*{\essinf}{ess\, inf}
\DeclareMathOperator*{\esssup}{ess\, sup}
\DeclareMathOperator*{\supp}{supp}

\DeclareMathOperator*{\a.e}{a.e.\; }

 \numberwithin{equation}{section}

 \usepackage[amsmath,thmmarks]{ntheorem} 

\theoremstyle{definition} 
\theoremheaderfont{\rmfamily\bf}
\theorembodyfont{\rm}

\newtheorem{theorem}{\indent
                  Theorem}[section]
    \newtheorem{lemma}{\indent  Lemma} [section]
\newtheorem{corollary} {\indent Corollary}[section]
\newtheorem{proposition}[theorem]{\indent  Proposition}

\theoremstyle{definition} 
    \newtheorem{definition}{\indent  Definition} [section]
        \newtheorem{remark}{\indent  Remark}  

\newtheorem*{acknowledgments}{\indent Acknowledgments}
\newtheorem*{proof}{\indent Proof}

\title{\bf\Large
 Boundedness of some sublinear operators  on  Herz-Morrey spaces with variable exponent  }

\author{
{\normalsize   Jianglong WU 
}
\\
{\small\it   Department of Mathematics, Mudanjiang Normal University, Mudanjiang \ 157011, China}\\
 {\small\it  Georgian Math. J. 2014; 21 (1):101-111}
  }
\date{} 

\begin{document}

\maketitle



\begin{minipage}[t]{14cm}

\setlength{\baselineskip}{1.0em}

\noindent

 { \bf Abstract:} In this paper, the boundedness of some sublinear operators is proved  on homogeneous Herz-Morrey spaces with variable exponent.


 \smallskip
 { \bf AMS(2010) Subject Classification:}  \ 42B20; 42B35   \ \

\end{minipage}

\section{Introduction}

Function spaces with variable exponent are being watched with keen interest not in real analysis but also in partial differential equations and in applied mathematics because they are applicable to the modeling for electrorheological fluids and image restoration. The theory of function spaces with variable exponent has rapidly made progress in the past twenty years
since some elementary properties were established by Kov\'{a}\v{c}ik-R\'{a}kosn\'{i}k\cite{KR}. One of the main
problems on the theory is the boundedness of the Hardy-Littlewood maximal operator on
variable Lebesgue spaces. By virtue of the fine works\cite{CDF,CFMP,CFN,D1,D2,DHHMS,K,L,N,PR}, some
important conditions on variable exponent, for example, the $\log$-H\"{o}lder conditions and the
Muckenhoupt type condition, have been obtained.

The class of the Herz spaces is arising from the study on characterization of multipliers on the classical Hardy spaces. The well-known Morrey spaces is used to show that certain systems of partial differential equations (PDEs) had H\"{o}lder
continuous solutions.  And the homogeneous Herz-Morrey spaces $M\dot{K}^{\alpha,\lambda}_{p,q}(\R^{n})$ coordinate with  the homogeneous Herz space $\dot{K}^{\alpha,p}_{q}(\R^{n})$ when $\lambda=0$.
One of the important problems on Herz spaces and Herz-Morrey spaces is the boundedness of sublinear
operators. Hern\'{a}ndez, Li, Lu and Yang\cite{HY,LiY,LY} have proved
that if a sublinear operator $T$ is bounded on $L^{p}(\R^{n})$ and satisfies the size condition
\begin{equation}  \label{Size}
 |Tf(x)|\le C \int_{\R^{n}} |x-y|^{-n}|f(y)|\mathrm{d}y,\ \ \ \a.e \ x\notin\supp f
 \end{equation}
for all $f\in L^{1}(\R^{n})$ with compact support, then $T$ is bounded on the homogeneous Herz space $\dot{K}^{\alpha,p}_{q}(\R^{n})$. In 2005, Lu and Xu\cite{LX} established the boundedness for some sublinear operators on homogeneous Herz-Morrey spaces.

In 2010, Izuki\cite{I2} proves the boundedness of some sublinear operators on Herz spaces with variable exponent; and recently Izuki\cite{I,I1} also considers the boundedness of some operators on Herz-Morrey spaces with variable exponent.

Motivated by the study on the Herz spaces and Lebesgue spaces  with variable exponent, when (\ref{Size}) is replaced by some more general size conditions, the main purpose of this paper is to establish some boundedness results of sublinear operators on Herz-Morrey spaces with variable exponent. This size condition is satisfied by many important operators in harmonic analysis.

Let us explain the outline of this article. In Section 2 we state some important properties of $ L^{q(\cdot)}(\R^{n})$ based on \cite{KR,D1,N,I}, and give some lemmas which will be needed for proving our main theorem. In Section 3 we prove the boundedness of
sublinear operators on Herz-Morrey spaces with variable exponent $ M\dot{K}_{p, q(\cdot)}^{\alpha, \lambda}(\mathbb{R}^{n} ) $, and obtain the corresponding corollaries.

Throughout this paper, we will denote by $|S|$ the Lebesgue measure  and  by $\chi_{_{\scriptstyle S}}$ the characteristic function  for a measurable set  $S\subset\R^{n}$. $C$  denotes a constant that is independent of the main parameters involved but
whose value may differ from line to line. For any index $1< q(x)< \infty$, we denote by $q'(x)$ its conjugate index,
namely, $q'(x)=\frac{q(x)}{q(x)-1}$.  For $A\sim D$, we mean that there is a constant $C > 0$ such that$C^{-1}D\le A \le CD$.

\section{
Preliminaries  and Lemmas}

In this section, we give the definition of Lebesgue and Herz-Morrey spaces with variable exponent, and
state their properties. Let $E$ be a measurable set in $\R^{n}$ with $|E|>0$. We first define Lebesgue spaces with variable exponent.

\begin{definition} \label{Def.L}
Let ~$ q(\cdot): E \to[1,\infty)$ be a measurable function.

1)\ The Lebesgue spaces with variable exponent $L^{q(\cdot)}(E)$ is defined by
  $$ L^{q(\cdot)}(E)=\{f~ \mbox{is measurable function}: \int_{E} \Big( \frac{|f(x)|}{\eta} \Big)^{q(x)} \mathrm{d}x <\infty ~\mbox{for some constant}~ \eta>0\}.  $$

  2)\ The space $L_{\loc}^{q(\cdot)}(E)$ is defined by
  $$ L_{\loc}^{q(\cdot)}(E)=\{f ~\mbox{is measurable function}: f\in L^{q(\cdot)}(K) ~\mbox{for all compact  subsets}~ K\subset E\}.  $$
\end{definition}

$L^{q(\cdot)}(E)$ is a Banach space with the norm defined by  $$ \|f\|_{L^{q(\cdot)}(E)}=\inf \Big\{ \eta>0:  \int_{E} \Big( \frac{|f(x)|}{\eta} \Big)^{q(x)} \mathrm{d}x \le 1 \Big\}.$$

Now, we define two classes of exponent functions. Given a function $f\in L_{\loc}^{1}(E)$, the Hardy-Littlewood maximal operator $M$ is defined by
$$Mf(x)=\sup_{r>0} r^{-n} \int_{B(x,r)\cap E} |f(y)| \mathrm{d}y \ \ \ (x\in E),$$
where $B(x,r)=\{y\in \R^{n}: |x-y|<r\}$.

\begin{definition} \label{Def.PB}
1)\ The set $\mathscr{P}(E)$ consists of all measurable functions $ q(\cdot)$ satisfying
$$1< \essinf_{x\in E} q(x)=q_{-},\ \ q_{+}= \esssup_{x\in E} q(x)<\infty.$$

2)\  The set $\mathscr{B}(E)$ consists of all  measurable functions  $q(\cdot)\in\mathscr{P}(E)$ such that the Hardy-Littlewood maximal operator $M$ is bounded on $L^{q(\cdot)}(E)$.
\end{definition}

Next we define the Herz-Morrey spaces with  variable exponent. Let $B_{k}=B(0,2^{k})=\{x\in\R^{n}:|x|\leq 2^{k}\}, A_{k}=\ B_{k}\setminus  B_{k-1}$ and $\chi_{_{k}}=\chi_{_{A_{k}}}$ for $k\in \mathbb{Z}$.

\begin{definition} \label{Def.M-H}
 Let $ \alpha\in \mathbb{R}, ~0\leq \lambda < \infty,~ 0<p< \infty$,
and $q(\cdot)\in \mathscr{P}(\mathbb{R}^{n})$. The Herz-Morrey space with  variable exponent $ M\dot{K}_{p, q(\cdot)}^{\alpha,
\lambda}(\mathbb{R}^{n} ) $ is defined by
  $$ M\dot{K}_{p, q(\cdot)}^{\alpha, \lambda}(\mathbb{R}^{n})=\{f\in L_{\loc}^{q(\cdot)}(\mathbb{R}^{n}\backslash\{0\}):
 \|f\|_{M\dot{K}_{p, q(\cdot)}^{\alpha, \lambda}(\mathbb{R}^{n})}<\infty\},  $$
where
  $ \|f\|_{M\dot{K}_{p, q(\cdot)}^{\alpha, \lambda}(\mathbb{R}^{n})}=\sup_{k_{0}\in \mathbb{Z}}2^{-k_{0}\lambda}
 \Big(\sum_{k=-\infty}^{k_{0}}2^{k\alpha
 p}\|f\chi_{_{\scriptstyle k}}\|_{L^{^{q(\cdot)}}(\mathbb{R}^{n})}^{p} \Big)^{\frac{1}{p}}.  $

\end{definition}

Compare the Herz-Morrey space with  variable exponent  $ M\dot{K}_{p, q(\cdot)}^{\alpha, \lambda}(\mathbb{R}^{n})$ with the Herz space with  variable exponent $\dot{K}_{q(\cdot)}^{\alpha,p}(\R^{n})$\cite{I1}, where
$$\dot{K}_{q(\cdot)}^{\alpha,p}(\R^{n})= \Big\{f\in L_{\loc}^{q(\cdot)}(\mathbb{R}^{n}\backslash\{0\}):
\sum\limits_{k=-\infty}^{\infty}2^{k\alpha  p}\|f\chi_{_{k}}\|_{L^{q(\cdot)}(\R^{n})}^{p}<\infty \Big\}.$$
Obviously, $M\dot{K}_{p,q(\cdot)}^{\alpha,0} (\R^{n})=\dot{K}_{q(\cdot)}^{\alpha,p}(\R^{n})$.

In 2012,  Almeida and  Drihem\cite{AD} discussed the boundedness of a wide class of sublinear operators, including maximal, potential and Calder\'{o}n-Zygmund operators, on variable Herz spaces $K_{q(\cdot)}^{\alpha(\cdot),p}(\R^{n})$ and $\dot{K}_{q(\cdot)}^{\alpha(\cdot),p}(\R^{n})$. Meanwhile, they also established Hardy-Littlewood-Sobolev theorems for fractional integrals on variable Herz spaces. In other papers \cite{AHS,KM1,KM2} the boundedness of operators of harmonic analysis are established in variable exponent Morrey spaces. In this paper, the author  considers Herz-Morrey space $ M\dot{K}_{p, q(\cdot)}^{\alpha(\cdot), \lambda}(\mathbb{R}^{n})$ with variable exponent $q(\cdot)$ for fixed $ \alpha\in \mathbb{R}$ and $p\in(0,\infty)$. However, for  the case of the exponent $\alpha(\cdot)$ is variable as well, we can refer to the furthermore work of the author of this paper.

Next we state some properties of variable exponent. Cruz-Uribe et al\cite{CFN} and Nekvinda\cite{N} proved the following sufficient conditions for the boundedness of $M$ in variable exponent space independently. We note that Nekvinda\cite{N} gave a more general condition in place of (\ref{C2}).

\begin{proposition}\hspace{-2pt}\cite{N}\quad   \label{Pr-1}
Suppose that $E$ is an open set. If $ q(\cdot)\in \mathscr{P}(E)$ satisfies the inequality
 \begin{alignat}{2}
 |q(x)-q(y)| &\le \frac{-C}{\ln(|x-y|)}  \ \ \ \mbox{if}~ |x-y|\le 1/2,   \label{C1}\\
 |q(x)-q(y)| &\le \frac{C}{\ln(e+|x|)}  \ \ \ \mbox{if}~ |y|\ge|x|,\ \    \label{C2}
\end{alignat}
where $C>0$ is a constant independent of $x$ and $y$, then we have $ q(\cdot)\in \mathscr{B}(E)$.
\end{proposition}

The next proposition is due to Diening \cite{D1}(see Theorem 8.1). We remark that Diening has also proved general results on Musielak-Orlicz spaces. We only describe partial results  we need in this paper.

\begin{proposition}\hspace{-2pt}\cite{D1}\quad \label{Pr-2}
Suppose that $ q(\cdot)\in \mathscr{P}(\R^{n})$, then
$q(\cdot)\in\mathscr{B}(\R^{n})$ iff $q'(\cdot)\in \mathscr{B}(\R^{n})$.


\end{proposition}

In order to prove our main theorems, we also need the following result which is the Hardy-Littlewood-Sobolev theorem on Lebesgue spaces with varible expoonent due to Capone, Cruz-Uribe and Fiorenza\cite{CCF}(see Theorem 1.8). We remark that this result was proved by Diening\cite{D} provided that $q_{_{1}}(\cdot)$ is constant outside of a large ball.

\begin{proposition}\hspace{-2pt}\cite{CCF} \quad   \label{Pr-3}
Suppose that $q_{_{1}}(\cdot)\in \mathscr{P}(\R^{n})$ satisfies conditions (\ref{C1}) and (\ref{C2}) in Proposition \ref{Pr-1}.
$0<\beta<n/(q_{_{1}})_{+}$ and define  $q_{_{2}}(\cdot)$ by
$$\frac{1}{q_{_{1}}(x)} -\frac{1}{q_{_{2}}(x)}=\frac{\beta}{n}.
$$
Then we have
$$\|I_{\beta}f\|_{L^{q_{2}(\cdot)}(\R^{n})} \le {C}\|f\|_{L^{q_{1}(\cdot)}(\R^{n})},$$
where  fractional integral $I_{\beta}$ is defined by
$I_{\beta}(f)(x) = \int_{\R^{n}}\frac{f(y)}{|x-y|^{n-\beta}} \mathrm{d}y$.
\end{proposition}

In  addition, by Proposition \ref{Pr-3}, we can obtain

\begin{proposition} \  \label{Pr-4}
Let $q_{_{1}}(\cdot),~ q_{_{2}}(\cdot)$ and $\beta$ be the same as in Proposition \ref{Pr-3}.
Then we have
$$\|\chi_{_{B_{k}}}\|_{L^{q_{2}(\cdot)}(\R^{n})} \le {C} 2^{-k\beta} \|\chi_{_{B_{k}}}\|_{L^{q_{1}(\cdot)}(\R^{n})}$$
for all balls $B_{k}=\{x\in\R^{n}:|x|\leq 2^{k}\}$ with  $k\in \mathbb{Z}$.
\end{proposition}

\begin{proof}
It is easy to see that
\begin{equation*}
\begin{split}
I_{\beta}(\chi_{B_{k}})(x) &\ge  I_{\beta}(\chi_{B_{k}})(x)\cdot \chi_{B_{k}}(x) = \int_{B_{k}}\frac{\mathrm{d}y }{|x-y|^{n-\beta}} \cdot \chi_{B_{k}}(x)\ge C 2^{k\beta}\chi_{B_{k}}(x).
\end{split}
\end{equation*}
Thus, applying Proposition \ref{Pr-3}, we get
\begin{eqnarray*}
\begin{split} \label{F-K}
\|\chi_{_{\scriptstyle B_{k}}}\|_{L^{^{q_{_{2}}(\cdot)}}(\R^{n})}      \le C 2^{^{-k\beta}} \|I_{\beta}(\chi_{_{\scriptstyle{B_{k}}}}) \|_{L^{^{q_{_{2}}(\cdot)}}(\R^{n})}  \le C2^{^{-k\beta}} \|\chi_{_{\scriptstyle B_{k}}} \|_{L^{^{q_{_{1}}(\cdot)}}(\R^{n})}.
 \end{split}
\end{eqnarray*}
This completes the proof of  Proposition \ref{Pr-4}.
\end{proof}

The next lemma describes the generalized H\"{o}lder's inequality and the
duality of $L^{q(\cdot)}(E)$. The proof can be found in \cite{KR}.

\begin{lemma} \hspace{-2pt}\cite{KR} \quad   \label{Lem.1}
Suppose that $ q(\cdot)\in \mathscr{P}(E)$. Then the following statements hold.

1)\  (generalized H\"{o}lder's inequality)\ \ For all $f\in L^{q(\cdot)}(E)$ and all $g\in L^{q'(\cdot)}(E)$, we have
 \begin{equation}  \label{Hd}
 \int_{E} |f(x)g(x)|\mathrm{d}x \le r_{q}\|f\|_{L^{q(\cdot)}(E)} \|g\|_{L^{q'(\cdot)}(E)},
\end{equation}
where $r_{q}=1+1/q_{-}-1/q_{+}$.

2)\ For all $f\in L^{q(\cdot)}(E)$, we have
$$\|f\|_{L^{q(\cdot)}(E)} \le \sup \Big\{ \int_{E}|f(x)g(x)|\mathrm{d}x: \|g\|_{L^{q'(\cdot)}(E)}\le 1 \Big\}.$$

\end{lemma}

\begin{lemma} \hspace{-2pt}\cite{I}\quad  \label{Lem.2}
If $ q(\cdot)\in \mathscr{B}(\R^{n})$, then there exist positive constants $\delta\in(0,1)$ and $C>0$ such that
 \begin{equation}  \label{BZ-1}
 \frac{\|\chi_{S}\|_{L^{q(\cdot)}(\R^{n})}}{\|\chi_{B}\|_{L^{q(\cdot)}(\R^{n})}} \le C \Big(\frac{|S|}{|B|} \Big)^{\delta}
\end{equation}
holds for all balls $B$ in $\R^{n}$ and all measurable subsets $S\subset B$.
\end{lemma}

\begin{lemma}\hspace{-2pt}\cite{I}\quad   \label{Lem.3}
If $ q(\cdot)\in \mathscr{B}(\R^{n})$, then there exists a positive constant  $C>0$ such that
 \begin{equation}  \label{equ.1}
 C^{-1}\le \frac{1}{|B|} \|\chi_{B}\|_{L^{q(\cdot)}(E)} \|\chi_{B}\|_{L^{q'(\cdot)}(E)} \le C.
\end{equation}
holds for all balls $B$ in $\R^{n}$.
\end{lemma}

\section{Main theorems and their proofs}

Let $ q(\cdot)\in \mathscr{P}(\R^{n})$ satisfy conditions (\ref{C1}) and (\ref{C2}) in Proposition \ref{Pr-1}. Then so does $ q'(\cdot)$. In particular, we can see that  $ q(\cdot),~q'(\cdot)\in \mathscr{B}(\R^{n})$ from  Proposition \ref{Pr-1}. Therefore applying Lemma \ref{Lem.2}, we can take  constants $\delta_{1}\in (0, 1/(q')_{+}), \delta_{2}\in (0, 1/(q)_{+})$ such that
\begin{equation}  \label{BZ}
 \frac{\|\chi_{S}\|_{L^{q'(\cdot)}(\R^{n})}} {\|\chi_{B}\|_{L^{q'(\cdot)}(\R^{n})}} \le C \bigg(\frac{|S|}{|B|} \bigg)^{\delta_{1}},
 \ \ \  \frac{\|\chi_{S}\|_{L^{q(\cdot)}(\R^{n})}} {\|\chi_{B}\|_{L^{q(\cdot)}(\R^{n})}} \le C \bigg(\frac{|S|}{|B|} \bigg)^{\delta_{2}}
\end{equation}
holds for all balls $B$ in $\R^{n}$ and all measurable subsets $S\subset B$. And when $q_{1}(\cdot), q_{2}(\cdot)\in \mathscr{P}(\R^{n})$, applying Lemma \ref{Lem.2}, we can take  constants $\delta_{3}\in (0, 1/(q'_{2})_{+}), \delta_{4}\in (0, 1/(q_{1})_{+})$ such that
\begin{equation}  \label{equ.4}
 \frac{\|\chi_{S}\|_{L^{q'_{1}(\cdot)}(\R^{n})}} {\|\chi_{B}\|_{L^{q'_{1}(\cdot)}(\R^{n})}} \le C \bigg(\frac{|S|}{|B|} \bigg)^{\delta_{3}},
 \ \ \  \frac{\|\chi_{S}\|_{L^{q_{2}(\cdot)}(\R^{n})}} {\|\chi_{B}\|_{L^{q_{2}(\cdot)}(\R^{n})}} \le C \bigg(\frac{|S|}{|B|} \bigg)^{\delta_{4}}
\end{equation}
holds for all balls $B$ in $\R^{n}$ and all measurable subsets $S\subset B$.

In this section, we will give some size condition which are more general than (\ref{Size}), and prove the boundedness of some sublinear operators, satisfying this size condition on Herz-Morrey spaces with variable exponent. 
This size condition is satisfied by many important operators in harmonic analysis.  Our main result can be stated as follows.

\begin{theorem}\label{thm.1}
Let $q(\cdot)\in \mathscr{P}(\R^{n})$ satisfies conditions (\ref{C1}) and (\ref{C2}) in Proposition \ref{Pr-1}, and $ 0 < p <\infty, ~ \lambda> 0,~\lambda-n\delta_{2}<\alpha<\lambda+ n\delta_{1}$,  where $\delta_{1}\in (0, 1/(q')_{+})$ and $\delta_{2}\in (0, 1/(q)_{+})$ are the constants satisfying (\ref{BZ}). Suppose that  a sublinear operator $T$ satisfies

\begin{list}{}{\setlength{\leftmargin}{3em}}
\item[(i)]  \ $T$ is bounded on $L^{q(\cdot)}(\R^{n})$;
\item[(ii)] \ for suitable function $f$ with $\supp f\subset A_{k}$ and $|x|\ge 2^{k+1}$ with $k\in \mathbb{Z}$,
\begin{equation}\label{Size-1}
|Tf(x)|\le C|x|^{-n}\|f\|_{L^{1}(\R^{n})};
\end{equation}
 \item[(iii)] \ for suitable function $f$ with $\supp f\subset A_{k}$ and $|x|\le 2^{k-2}$ with $k\in \mathbb{Z}$,
\begin{equation}\label{Size-2}
|Tf(x)|\le C2^{-kn}\|f\|_{L^{1}(\R^{n})}.
\end{equation}
\end{list}
Then $T$ is also bounded on $M\dot{K}_{p,q(\cdot)}^{\alpha,\lambda}(\R^{n})$.

\end{theorem}

Note that if $T$ satisfies the size condition (\ref{Size}), then $T$ satisfies (\ref{Size-1}) and (\ref{Size-2}). Thus, by Theorem \ref{thm.1}, we have
\begin{corollary}\label{co.1}
Let $q(\cdot),~p,~\lambda,~\alpha,~\delta_{1}$ and $\delta_{2}$ be the same as in Theorem \ref{thm.1}. If a sublinear operator $T$ satisfies the size condition (\ref{Size}) and is bounded on $L^{q(\cdot)}(\R^{n})$, then  $T$ is also bounded on $M\dot{K}_{p,q(\cdot)}^{\alpha,\lambda}(\R^{n})$.

\end{corollary}

\begin{remark}
For any $q(\cdot)\in \mathscr{P}(\R^{n})$ satisfies conditions (\ref{C1}) and (\ref{C2}) in Proposition \ref{Pr-1}, the Hardy-Littlewood maximal operator $M$, defined by
$$M(f)(x)=\sup_{B\ni x}\frac{1}{|B|} \int_{B\cap \Omega}|f(y)|\mathrm{d}y, $$
also satisfies the assumptions  of Theorem \ref{thm.1}, where the supremum is  taken over all balls $B$ containing $x$, and $\Omega\subset \R^{n}$ is an open set.
\end{remark}

 \begin{proof}[\bf{\it Proof of Theorem \ref{thm.1}}]
For all $  f\in {M\dot{K}_{p,  q(\cdot)}^{\alpha,\lambda}(\R^{n})}$. If we denote   $f_j:=f\cdot\chi_{j}=f\cdot\chi_{A_j}$ for each $j\in \mathbb{Z}$, then we can write
$$f(x)=\sum_{j=-\infty}^{\infty}f(x)\cdot\chi_{j}(x) =\sum_{j=-\infty}^{\infty}f_{j}(x).
$$
We have
\begin{eqnarray*}
&&\; \|Tf\|^{p}_{M\dot{K}_{p, q(\cdot)}^{\alpha,\lambda}(\R^{n})}= \sup_{k_{0}\in \mathbb{Z}}2^{-k_{0}\lambda p}  \bigg(\sum_{k=-\infty}^{k_{0}}2^{k\alpha p}  \| T(f)\cdot\chi_{_{\scriptstyle k}} \|_{L^{q(\cdot)}(\R^{n})}^{p}\bigg)\\
&&\; \le C\sup_{k_{0}\in \mathbb{Z}}2^{-k_{0}\lambda p}  \bigg(\sum_{k=-\infty}^{k_{0}}2^{k\alpha p}   \Big(\sum_{j=-\infty}^{k-2}\| T(f_{j})\cdot\chi_{_{\scriptstyle k}} \|_{L^{q(\cdot)}(\R^{n})}\Big)^{p}\bigg)\\
&&\;\ \ + C\sup_{k_{0}\in \mathbb{Z}}2^{-k_{0}\lambda p}  \bigg(\sum_{k=-\infty}^{k_{0}}2^{k\alpha p}   \Big\| T \Big(\sum_{j=k-1}^{k+1}f_{j}\Big)\cdot\chi_{_{\scriptstyle k}} \Big\|_{L^{q(\cdot)}(\R^{n})}^{p}\bigg)\\
&&\;\ \ + C\sup_{k_{0}\in \mathbb{Z}}2^{-k_{0}\lambda p}  \bigg(\sum_{k=-\infty}^{k_{0}}2^{k\alpha p}   \Big(\sum_{j=k+2}^{\infty}\| T(f_{j})\cdot\chi_{_{\scriptstyle k}} \|_{L^{q(\cdot)}(\R^{n})}\Big)^{p}\bigg)\\
 &&\;  = C(E_{1} +E_{2}+E_{3}).
\end{eqnarray*}

First we estimate $E_{2}$. Applying the $L^{q(\cdot)}(\R^{n})$-boundedness of $T$, we have
\begin{eqnarray*}
 E_{2} &=&  \sup_{k_{0}\in \mathbb{Z}}2^{-k_{0}\lambda p}  \bigg(\sum_{k=-\infty}^{k_{0}}2^{k\alpha p}   \Big\| T \Big(\sum_{j=k-1}^{k+1}f_{j}\Big)\cdot\chi_{_{\scriptstyle k}} \Big\|_{L^{q(\cdot)}(\R^{n})}^{p}\bigg)\\
&\le&  C\sup_{k_{0}\in \mathbb{Z}}2^{-k_{0}\lambda p}  \bigg(\sum_{k=-\infty}^{k_{0}}2^{k\alpha p}   \Big\| \Big(\sum_{j=k-1}^{k+1}f_{j}\Big)\cdot\chi_{_{\scriptstyle k}} \Big\|_{L^{q(\cdot)}(\R^{n})}^{p}\bigg)\\
&\le& C\sup_{k_{0}\in \mathbb{Z}}2^{-k_{0}\lambda p}  \bigg(\sum_{k=-\infty}^{k_{0}}2^{k\alpha p}   \| f\cdot\chi_{_{\scriptstyle k}} \|_{L^{q(\cdot)}(\R^{n})}^{p}\bigg)\\
 &=& C\|f\|^{p}_{M\dot{K}_{p,q(\cdot)}^{\alpha,\lambda}(\R^{n})}.
\end{eqnarray*}

For $E_{1}$, we  notice the facts that  $ j\le k-2$ and $\a.e x\in A_{k}$ with $k\in \mathbb{Z}$, then using the size condition (\ref{Size-1}) and the generalized H\"{o}lder's inequality(see (\ref{Hd}) in Lemma \ref{Lem.1}), we have
\begin{equation}\label{Size-11}
\begin{split}
|T(f_{j})(x)|&\le C|x|^{-n}\|f_{j}\|_{L^{1}(\R^{n})}\le C2^{-kn}\|f_{j}\|_{L^{1}(\R^{n})}\\
 &\le C 2^{-kn}\|f_{j}\|_{L^{q(\cdot)}(\R^{n})}\| \chi_{_{\scriptstyle j}}\|_{L^{q'(\cdot)}(\R^{n})}.
\end{split}
\end{equation}
Using  Proposition \ref{Pr-1}, Proposition \ref{Pr-2}, Lemma \ref{Lem.2},  Lemma \ref{Lem.3} and (\ref{BZ}), we obtain
\begin{equation}\label{equ.2}
\begin{split}
2^{-kn}\|\chi_{_{\scriptstyle k}}\|_{L^{q(\cdot)}(\R^{n})} \|\chi_{_{\scriptstyle j}}\|_{L^{q'(\cdot)}(\R^{n})}
&\le  2^{-kn}\|\chi_{_{B_{k}}}\|_{L^{q(\cdot)}(\R^{n})} \|\chi_{_{B_{j}}}\|_{L^{q'(\cdot)}(\R^{n})} \\
&\le C \|\chi_{_{B_{k}}}\|^{-1}_{L^{q'(\cdot)}(\R^{n})} \|\chi_{_{B_{j}}}\|_{L^{q'(\cdot)}(\R^{n})} \\
 &= C \frac{\|\chi_{_{B_{j}}}\|_{L^{q'(\cdot)}(\R^{n})}}{\|\chi_{_{B_{k}}}\|_{L^{q'(\cdot)}(\R^{n})}}\le C 2^{(j-k)n\delta_{1}}.
\end{split}
\end{equation}

On the other hand, note the following fact
\begin{eqnarray}
\begin{split} \label{Fj}
\|f_{j}\|_{L^{^{q(\cdot)}}(\R^{n})} &= 2^{-j\alpha}\Big(2^{j{\alpha}p} \|f_{j}\|^{p}_{L^{^{q(\cdot)}}(\R^{n})}\Big)^{1/p}\\
&\le  2^{-j\alpha}\bigg(\sum_{i=-\infty}^j 2^{i\alpha p} \|f_{i}\|^{p}_{L^{^{q}(\cdot)}(\R^{n})}\bigg)^{1/p}\\
&=  2^{j(\lambda-\alpha)}\bigg(2^{-j\lambda} \Big(\sum_{i=-\infty}^j 2^{i{\alpha}p} \|f_{i}\|^{p}_{L^{^{q(\cdot)}} (\R^{n})}\Big)^{1/p}\bigg)\\
&\le  C 2^{j(\lambda-\alpha)} \|f\|_{M\dot{K}_{p,q(\cdot)}^{\alpha,\lambda}(\R^{n})}.
 \end{split}
\end{eqnarray}
 Therefore, combining (\ref{Size-11}), (\ref{equ.2}) and (\ref{Fj}), and using $\alpha< \lambda+n\delta_{1}$, it follows that
\begin{equation*}
\begin{split}
E_{1}  & \le C\sup_{k_{0}\in \mathbb{Z}}2^{-k_{0}\lambda p}  \bigg(\sum_{k=-\infty}^{k_{0}}2^{k\alpha p}   \Big(\sum_{j=-\infty}^{k-2} 2^{-kn}\|f_{j}\|_{L^{q(\cdot)}(\R^{n})}\| \chi_{_{\scriptstyle j}}\|_{L^{q'(\cdot)}(\R^{n})} \| \chi_{_{\scriptstyle k}}\|_{L^{q(\cdot)}(\R^{n})}\Big)^{p}\bigg)  \\
 & \le C\sup_{k_{0}\in \mathbb{Z}}2^{-k_{0}\lambda p}  \bigg(\sum_{k=-\infty}^{k_{0}}2^{k\alpha p}   \Big(\sum_{j=-\infty}^{k-2} 2^{(j-k)n\delta_{1}}\|f_{j}\|_{L^{q(\cdot)}(\R^{n})} \Big)^{p}\bigg)  \\
 & \le C\sup_{k_{0}\in \mathbb{Z}}2^{-k_{0}\lambda p}  \bigg(\sum_{k=-\infty}^{k_{0}}2^{k\alpha p}   \Big(\sum_{j=-\infty}^{k-2} 2^{(j-k)n\delta_{1}}2^{j(\lambda-\alpha)} \|f\|_{M\dot{K}_{p,q(\cdot)}^{\alpha,\lambda}(\R^{n})} \Big)^{p}\bigg)  \\
 & \le C\|f\|^{p}_{M\dot{K}_{p,q(\cdot)}^{\alpha,\lambda}(\R^{n})}\sup_{k_{0}\in \mathbb{Z}}2^{-k_{0}\lambda p}  \bigg(\sum_{k=-\infty}^{k_{0}}2^{k\lambda p}   \Big(\sum_{j=-\infty}^{k-2} 2^{(j-k)(n\delta_{1}+\lambda-\alpha)} \Big)^{p}\bigg)  \\
    & \le C\|f\|^{p}_{M\dot{K}_{p,q(\cdot)}^{\alpha,\lambda}(\R^{n})}\sup_{k_{0}\in \mathbb{Z}}2^{-k_{0}\lambda p}  \bigg(\sum_{k=-\infty}^{k_{0}}2^{k\lambda p}   \bigg)  \le C\|f\|^{p}_{M\dot{K}_{p,q(\cdot)}^{\alpha,\lambda}(\R^{n})}.
\end{split}
\end{equation*}

Now, let us   estimate   $E_{3}$. For every $ j\ge k+2$ and $\a.e x\in A_{k}$ with $k\in \mathbb{Z}$, applying the size condition (\ref{Size-2}) and the generalized H\"{o}lder's inequality(see (\ref{Hd}) in Lemma \ref{Lem.1}), we have
\begin{equation}\label{Size-21}
\begin{split}
|T(f_{j})(x)|&\le  C2^{-jn}\|f_{j}\|_{L^{1}(\R^{n})} \le C 2^{-jn}\|f_{j}\|_{L^{q(\cdot)}(\R^{n})}\| \chi_{_{\scriptstyle j}}\|_{L^{q'(\cdot)}(\R^{n})}.
\end{split}
\end{equation}
Using Proposition \ref{Pr-1}, Lemma \ref{Lem.2}, Lemma \ref{Lem.3} and (\ref{BZ}), we obtain
\begin{equation}\label{equ.3}
\begin{split}
2^{-jn}\|\chi_{_{\scriptstyle k}}\|_{L^{q(\cdot)}(\R^{n})} \|\chi_{_{\scriptstyle j}}\|_{L^{q'(\cdot)}(\R^{n})}
&\le  2^{-jn}\|\chi_{_{B_{k}}}\|_{L^{q(\cdot)}(\R^{n})} \|\chi_{_{B_{j}}}\|_{L^{q'(\cdot)}(\R^{n})} \\
&\le C \|\chi_{_{B_{k}}}\|_{L^{q(\cdot)}(\R^{n})} \|\chi_{_{B_{j}}}\|^{-1}_{L^{q(\cdot)}(\R^{n})} \\
 &= C \frac{\|\chi_{_{B_{k}}}\|_{L^{q(\cdot)}(\R^{n})}}{\|\chi_{_{B_{j}}}\|_{L^{q(\cdot)}(\R^{n})}}\le C 2^{(k-j)n\delta_{2}}.
\end{split}
\end{equation}
Thus,  combining (\ref{Fj}), (\ref{Size-21}) and (\ref{equ.3}), and using $\alpha> \lambda-n\delta_{2}$, it follows that
\begin{equation*}
\begin{split}
E_{3}&=\sup_{k_{0}\in \mathbb{Z}}2^{-k_{0}\lambda p}  \bigg(\sum_{k=-\infty}^{k_{0}}2^{k\alpha p}   \Big(\sum_{j=k+2}^{\infty}\| T(f_{j})\cdot\chi_{_{\scriptstyle k}} \|_{L^{q(\cdot)}(\R^{n})}\Big)^{p}\bigg) \\
  & \le C\sup_{k_{0}\in \mathbb{Z}}2^{-k_{0}\lambda p}  \bigg(\sum_{k=-\infty}^{k_{0}}2^{k\alpha p}   \Big(\sum_{j=k+2}^{\infty} 2^{-jn}\|f_{j}\|_{L^{q(\cdot)}(\R^{n})}\| \chi_{_{\scriptstyle j}}\|_{L^{q'(\cdot)}(\R^{n})} \| \chi_{_{\scriptstyle k}}\|_{L^{q(\cdot)}(\R^{n})}\Big)^{p}\bigg)  \\
 & \le C\sup_{k_{0}\in \mathbb{Z}}2^{-k_{0}\lambda p}  \bigg(\sum_{k=-\infty}^{k_{0}}2^{k\alpha p}   \Big(\sum_{j=k+2}^{\infty} 2^{(k-j)n\delta_{2}}\|f_{j}\|_{L^{q(\cdot)}(\R^{n})} \Big)^{p}\bigg)  \\
 & \le C\sup_{k_{0}\in \mathbb{Z}}2^{-k_{0}\lambda p}  \bigg(\sum_{k=-\infty}^{k_{0}}2^{k\alpha p}   \Big(\sum_{j=k+2}^{\infty} 2^{(k-j)n\delta_{2}} 2^{j(\lambda-\alpha)} \|f\|_{M\dot{K}_{p,q(\cdot)}^{\alpha,\lambda}(\R^{n})} \Big)^{p}\bigg)  \\
   & \le C\|f\|^{p}_{M\dot{K}_{p,q(\cdot)}^{\alpha,\lambda}(\R^{n})}\sup_{k_{0}\in \mathbb{Z}}2^{-k_{0}\lambda p}  \bigg(\sum_{k=-\infty}^{k_{0}}2^{k\lambda p}   \Big(\sum_{j=k+2}^{\infty} 2^{(j-k)(\lambda-\alpha-n\delta_{2})} \Big)^{p}\bigg)  \\
    & \le C\|f\|^{p}_{M\dot{K}_{p,q(\cdot)}^{\alpha,\lambda}(\R^{n})}\sup_{k_{0}\in \mathbb{Z}}2^{-k_{0}\lambda p}  \bigg(\sum_{k=-\infty}^{k_{0}}2^{k\lambda p}   \bigg)  \le C\|f\|^{p}_{M\dot{K}_{p,q(\cdot)}^{\alpha,\lambda}(\R^{n})}.
\end{split}      
\end{equation*}

Combining the estimates for  $E_{1}$,  $E_{2}$ and  $E_{3}$ yields
$$\|Tf\|_{M\dot{K}_{p, q(\cdot)}^{\alpha,\lambda}(\R^{n})} \le C\|f\|_{M\dot{K}_{p,q(\cdot)}^{\alpha,\lambda}(\R^{n})}$$
and then completes the proof of Theorem \ref{thm.1}.
\end{proof}

Now, let us turn to consider the fractional singular integrals. We first have the  following theorem similar to \ref{thm.1}.

\begin{theorem}\label{thm.2}
Let $q_{_{1}}(\cdot)\in \mathscr{P}(\R^{n})$ satisfies conditions (\ref{C1}) and (\ref{C2}) in Proposition \ref{Pr-1}. Define the variable exponent $q_{_{2}}(\cdot)$ by
$$\frac{1}{q_{_{1}}(x)} -\frac{1}{q_{_{2}}(x)} = \frac{\beta}{n}.$$
 And let $ 0<p_{_{1}}\le {p_{_{2}}} <\infty,~ \lambda> 0,~ 0<\beta<n/(q_{_{1}})_{+},~ \lambda-n\delta_{4}<\alpha<\lambda+ n\delta_{3}$,  where $\delta_{3}\in (0, 1/(q'_{1})_{+})$ and $\delta_{4}\in (0, 1/(q_{2})_{+})$ are the constants appearing in (\ref{equ.4}). Suppose that  a sublinear operator $T_{\beta}$ satisfies
\begin{list}{}{\setlength{\leftmargin}{3em}}
\item[(i)]  \ $T_{\beta}$ maps from $L^{^{q_{1}(\cdot)}}(\R^{n})$ to $L^{^{q_{2}(\cdot)}}(\R^{n})$;
\item[(ii)] \ for any function $f$ with $\supp f\subset A_{k}$ and any $|x|\ge 2^{k+1}$ with $k\in \mathbb{Z}$,
\begin{equation}\label{equ.5}
|T_{\beta}(f)(x)|\le C|x|^{\beta-n}\|f\|_{L^{1}(\R^{n})};
\end{equation}
 \item[(iii)] \ for any function $f$ with $\supp f\subset A_{k}$ and any $|x|\le 2^{k-2}$ with $k\in \mathbb{Z}$,
\begin{equation}\label{equ.6}
|T_{\beta}(f)(x)|\le C2^{k(\beta-n)}\|f\|_{L^{1}(\R^{n})}.
\end{equation}
\end{list}
Then  $T_{\beta}(f)$ is  bounded from $M\dot{K}_{p_{_{1}},q_{_{1}}(\cdot)}^{\alpha,\lambda}(\R^{n})$ to $M\dot{K}_{p_{_{2}},q_{_{2}}(\cdot)}^{\alpha,\lambda}(\R^{n})$.
\end{theorem}

Notice that if $T_{\beta}(f)$ satisfies the size condition
 \begin{equation}  \label{equ.12}
 |T_{\beta}(f)(x)|\le C \int_{\R^{n}}\frac{|f(y)|}{|x-y|^{n-\beta}}\mathrm{d}y, \ \ \ \ x\notin \supp f
 \end{equation}
for any integrable function $f$ with compact support, then $T_{\beta}(f)$ obviously satisfies the assumptions of Theorem \ref{thm.2}.  Therefore, by Theorem \ref{thm.2}, we have
\begin{corollary}\label{co.2}
Let $q_{_{1}}(\cdot),~q_{_{2}}(\cdot),~p_{_{1}},~ {p_{_{2}}},~\lambda,~\beta,~\alpha,~\delta_{3}$ and $\delta_{4}$ be the same as in Theorem \ref{thm.2}. If a sublinear operator $T_{\beta}$ satisfies the size condition (\ref{equ.12}) and is bounded from  $L^{^{q_{1}}}(\R^{n})$ to $L^{^{q_{2}}}(\R^{n})$, then  $T_{\beta}$ is also bounded from $M\dot{K}_{p_{_{1}},q_{_{1}}(\cdot)}^{\alpha,\lambda}(\R^{n})$ to $M\dot{K}_{p_{_{2}},q_{_{2}}(\cdot)}^{\alpha,\lambda}(\R^{n})$.
\end{corollary}

\begin{remark}
We can see that when $\beta=0$, Theorem \ref{thm.2} is just Theorem \ref{thm.1}.
In particular, if $T_{\beta}(f)$ is a (standard) fractional integral, then $T_{\beta}(f)$ obviously satisfies (\ref{equ.12}). The fractional maximal function $M_{\beta}(f)$, defined by
$$M_{\beta}(f)(x)=\sup_{B\ni x}\frac{1}{|B|^{1-\beta/n}} \int_{B\cap \Omega}|f(y)|\mathrm{d}y, $$
also satisfies the conditions of Theorem \ref{thm.2}, where the supremum is again taken over all balls $B$ which contain $x$,  and $\Omega\subset \R^{n}$ is an open set.
\end{remark}

\begin{proof}[\bf{\it Proof of Theorem \ref{thm.2}}]
For all $f\in {M\dot{K}_{p_{_{1}},  q_{_{1}}(\cdot)}^{\alpha,\lambda}(\R^{n})}$. If we denote $f_j:=f\cdot\chi_{j}=f\cdot\chi_{A_j}$ for each $j\in \mathbb{Z}$, then we can write
$$f(x)=\sum_{j=-\infty}^{\infty}f(x)\cdot\chi_{j}(x) =\sum_{j=-\infty}^{\infty}f_{j}(x).
$$
Because of $0<p_{_{1}}/p_{_{2}}\le 1$, we apply inequality
\begin{equation}     \label{Cp}
\bigg(\sum_{i=-\infty}^{\infty}|a_{i}|\bigg)^{ p_{_{1}}/p_{_{2}}} \le \sum_{i=-\infty}^{\infty} |a_{i}|^{ p_{_{1}}/ p_{_{2}}},
\end{equation}
and obtain
\begin{eqnarray*}
&&\; \|T_{\beta}(f)\|^{p_{_{1}}}_{M\dot{K}_{p_{_{2}}, q_{_{2}}(\cdot)}^{\alpha,\lambda}(\R^{n})}= \sup_{k_{0}\in \mathbb{Z}}2^{-k_{0}\lambda {p_{_{1}}}}  \bigg(\sum_{k=-\infty}^{k_{0}}2^{k\alpha {p_{_{2}}}}  \| T_{\beta}(f)\cdot\chi_{_{\scriptstyle k}}  \|_{L^{q_{_{2}}(\cdot)}(\R^{n})}^{p_{_{2}}}\bigg)^{p_{_{1}}/p_{_{2}}}\\
&&\; \le C \sup_{k_{0}\in \mathbb{Z}}2^{-k_{0}\lambda p_{_{1}}}  \bigg(\sum_{k=-\infty}^{k_{0}}2^{k\alpha p_{_{1}}}
  \|T_{\beta}(f)\cdot\chi_{_{\scriptstyle{k}}}  \|_{L^{q_{_{2}}(\cdot)}(\R^{n})}^{p_{_{1}}}\bigg) \\
&&\; \le C \sup_{k_{0}\in \mathbb{Z}}2^{-k_{0}\lambda p_{_{1}}}  \bigg(\sum_{k=-\infty}^{k_{0}}2^{k\alpha p_{_{1}}}
  \Big(\sum_{j=-\infty}^{k-2}\|T_{\beta}(f_{j})\cdot\chi_{_{\scriptstyle{k}}}  \|_{L^{q_{_{2}}(\cdot)}(\R^{n})} \Big)^{p_{_{1}}}\bigg) \\
&&\;\ \ + C \sup_{k_{0}\in \mathbb{Z}}2^{-k_{0}\lambda p_{_{1}}}  \bigg(\sum_{k=-\infty}^{k_{0}}2^{k\alpha p_{_{1}}}
  \Big\|T_{\beta}\Big(\sum_{j=k-1}^{k+1}f_{j}\Big)\cdot\chi_{_{\scriptstyle{k}}} \Big \|_{L^{q_{_{2}}(\cdot)}(\R^{n})} ^{p_{_{1}}}\bigg)\\
&&\;\ \ + C \sup_{k_{0}\in \mathbb{Z}}2^{-k_{0}\lambda p_{_{1}}}  \bigg(\sum_{k=-\infty}^{k_{0}}2^{k\alpha p_{_{1}}}
  \Big(\sum_{j=k+2}^{\infty}\|T_{\beta}(f_{j})\cdot\chi_{_{\scriptstyle{k}}}  \|_{L^{q_{_{2}}(\cdot)}(\R^{n})} \Big)^{p_{_{1}}}\bigg)\\
 &&\;  = C(E_{1} +E_{2}+E_{3}).
\end{eqnarray*}

For $E_{2}$, using the boundedness of $T_{\beta}$  from $L^{^{q_{1}(\cdot)}}(\R^{n})$ to $L^{^{q_{2}(\cdot)}}(\R^{n})$,
we have
\begin{eqnarray*}
&&\; E_{2} = \sup_{k_{0}\in \mathbb{Z}}2^{-k_{0}\lambda p_{_{1}}}  \bigg(\sum_{k=-\infty}^{k_{0}}2^{k\alpha p_{_{1}}}   \Big\|T_{\beta}\Big(\sum_{j=k-1}^{k+1}f_{j}\Big)\cdot\chi_{_{\scriptstyle{k}}} \Big \|_{L^{q_{_{2}}(\cdot)}(\R^{n})} ^{p_{_{1}}}\bigg)\\
&&\; \le C \sup_{k_{0}\in \mathbb{Z}}2^{-k_{0}\lambda p_{_{1}}}  \bigg(\sum_{k=-\infty}^{k_{0}}2^{k\alpha p_{_{1}}}
  \Big\|\Big(\sum_{j=k-1}^{k+1}f_{j}\Big)\cdot\chi_{_{\scriptstyle{k}}}  \Big\|_{L^{q_{_{1}}(\cdot)}(\R^{n})}^{p_{_{1}}}\bigg)\\
  &&\;\le C \sup_{k_{0}\in \mathbb{Z}}2^{-k_{0}\lambda p_{_{1}}}  \bigg(\sum_{k=-\infty}^{k_{0}}2^{k\alpha p_{_{1}}}\|f\cdot\chi_{_{\scriptstyle{k}}}  \|_{L^{q_{_{1}}(\cdot)}(\R^{n})} ^{p_{_{1}}}\bigg)\\
 &&\;  = C \|f\|^{p_{1}}_{M\dot{K}_{p_{_{1}},q_{_{1}}(\cdot)}^{\alpha,\lambda}(\R^{n})}.
\end{eqnarray*}

For $E_{1}$. Note that  $j\le k-2$ and $\a.e x\in A_{k}$ with $k\in \mathbb{Z}$, then using the size condition (\ref{equ.5}) and the generalized H\"{o}lder's inequality(see (\ref{Hd}) in Lemma \ref{Lem.1}), we have
\begin{equation}\label{equ.7}
\begin{split}
|T_{\beta}(f_{j})(x) |&\le C|x|^{\beta-n}\|f_{j}\|_{L^{1}(\R^{n})}\le C2^{k(\beta-n)}\|f_{j}\|_{L^{1}(\R^{n})}\\
 &\le C 2^{k(\beta-n)}\|f_{j}\|_{L^{q_{1}(\cdot)}(\R^{n})}\| \chi_{_{\scriptstyle j}}\|_{L^{q'_{1}(\cdot)}(\R^{n})}.
\end{split}
\end{equation}
Using Proposition \ref{Pr-1}, Proposition \ref{Pr-2}, Proposition \ref{Pr-4}, Lemma \ref{Lem.2}, Lemma \ref{Lem.3} and (\ref{equ.4}), we obtain
\begin{equation}\label{equ.8}
\begin{split}
 2^{k(\beta-n)}\|\chi_{_{\scriptstyle k}}\|_{L^{q_{2}(\cdot)}(\R^{n})} \|\chi_{_{\scriptstyle  j}}\|_{L^{q'_{1}(\cdot)}(\R^{n})}
 &\le C 2^{-kn}2^{k\beta}\|\chi_{_{B_{k}}}\|_{L^{q_{2}(\cdot)}(\R^{n})} \|\chi_{_{B_{j}}}\|_{L^{q'_{1}(\cdot)}(\R^{n})} \\
 &\le C 2^{-kn}\|\chi_{_{B_{k}}}\|_{L^{q_{1}(\cdot)}(\R^{n})} \|\chi_{_{B_{j}}}\|_{L^{q'_{1}(\cdot)}(\R^{n})} \\
&\le C \|\chi_{_{B_{k}}}\|^{-1}_{L^{q'_{1}(\cdot)}(\R^{n})} \|\chi_{_{B_{j}}}\|_{L^{q'_{1}(\cdot)}(\R^{n})} \\
 &= C \frac{\|\chi_{_{B_{j}}}\|_{L^{q'_{1}(\cdot)}(\R^{n})}}{\|\chi_{_{B_{k}}}\|_{L^{q'_{1}(\cdot)}(\R^{n})}}\le C 2^{(j-k)n\delta_{3}}.
\end{split}
\end{equation}
On the other hand, note the following fact
\begin{eqnarray}
\begin{split} \label{equ.9}
\|f_{j}\|_{L^{^{q_{1}(\cdot)}}(\R^{n})} &= 2^{-j\alpha}\Big(2^{j{\alpha}p_{1}} \|f_{j}\|^{p_{1}}_{L^{^{q_{1}(\cdot)}}(\R^{n})}\Big)^{1/p_{1}}\\
&\le  2^{-j\alpha}\bigg(\sum_{i=-\infty}^j 2^{i\alpha p_{1}} \|f_{i}\|^{p_{1}}_{L^{^{q_{1}}(\cdot)}(\R^{n})}\bigg)^{1/p_{1}}\\
&=  2^{j(\lambda-\alpha)}\bigg(2^{-j\lambda} \Big(\sum_{i=-\infty}^j 2^{i{\alpha}p_{1}} \|f_{i}\|^{p_{1}}_{L^{^{q_{1}(\cdot)}} (\R^{n})}\Big)^{1/p_{1}}\bigg)\\
&\le  C 2^{j(\lambda-\alpha)} \|f\|_{M\dot{K}_{p_{1},q_{1}(\cdot)}^{\alpha,\lambda}(\R^{n})}.
 \end{split}
\end{eqnarray}
Hence, combining (\ref{equ.7}), (\ref{equ.8}) and (\ref{equ.9}), and using $\alpha< \lambda+n\delta_{3}$, it follows that
\begin{equation*}
\begin{split}
E_{1}&=\sup_{k_{0}\in \mathbb{Z}}2^{-k_{0}\lambda p_{_{1}}}  \bigg(\sum_{k=-\infty}^{k_{0}}2^{k\alpha p_{_{1}}}\Big( \sum_{j=-\infty}^{k-2} \|T_{\beta}(f_{j})\cdot\chi_{_{\scriptstyle{k}}}  \|_{L^{q_{_{2}}(\cdot)}(\R^{n})} \Big)^{p_{_{1}}}\bigg)  \\
&\le C \sup_{k_{0}\in \mathbb{Z}}2^{-k_{0}\lambda p_{_{1}}}  \bigg(\sum_{k=-\infty}^{k_{0}}2^{k\alpha p_{_{1}}}\Big( \sum_{j=-\infty}^{k-2} 2^{k(\beta-n)}\|f_{j}\|_{L^{q_{1}(\cdot)}(\R^{n})}\| \chi_{_{\scriptstyle j}}\|_{L^{q'_{1}(\cdot)}(\R^{n})}\|\chi_{_{\scriptstyle{k}}}  \|_{L^{q_{_{2}}(\cdot)}(\R^{n})} \Big)^{p_{_{1}}}\bigg)  \\
&\le C \sup_{k_{0}\in \mathbb{Z}}2^{-k_{0}\lambda p_{_{1}}}  \bigg(\sum_{k=-\infty}^{k_{0}}2^{k\alpha p_{_{1}}}\Big( \sum_{j=-\infty}^{k-2} 2^{(j-k)n\delta_{3}}\|f_{j}\|_{L^{q_{1}(\cdot)}(\R^{n})}\Big)^{p_{_{1}}}\bigg)  \\
&\le C\|f\|^{p_{1}}_{M\dot{K}_{p_{1},q_{1}(\cdot)}^{\alpha,\lambda}(\R^{n})}\sup_{k_{0}\in \mathbb{Z}}2^{-k_{0}\lambda p_{_{1}}}  \bigg(\sum_{k=-\infty}^{k_{0}}2^{k\alpha p_{_{1}}}\Big( \sum_{j=-\infty}^{k-2} 2^{(j-k)n\delta_{3}}2^{j(\lambda-\alpha)}\Big)^{p_{_{1}}}\bigg)  \\
&\le C\|f\|^{p_{1}}_{M\dot{K}_{p_{1},q_{1}(\cdot)}^{\alpha,\lambda}(\R^{n})}\sup_{k_{0}\in \mathbb{Z}}2^{-k_{0}\lambda p_{_{1}}}  \bigg(\sum_{k=-\infty}^{k_{0}}2^{k\lambda p_{_{1}}}\Big( \sum_{j=-\infty}^{k-2} 2^{(j-k)(\lambda-\alpha+n\delta_{3})}\Big)^{p_{_{1}}}\bigg)  \\
&\le C\|f\|^{p_{1}}_{M\dot{K}_{p_{1},q_{1}(\cdot)}^{\alpha,\lambda}(\R^{n})}\sup_{k_{0}\in \mathbb{Z}}2^{-k_{0}\lambda p_{_{1}}}  \Big(\sum_{k=-\infty}^{k_{0}}2^{k\lambda p_{_{1}}} \Big)  \le C\|f\|^{p_{1}}_{M\dot{K}_{p_{1},q_{1}(\cdot)}^{\alpha,\lambda}(\R^{n})}.
\end{split}
\end{equation*}

Now, let us  estimate   $E_{3}$. For every $ j\ge k+2$ and $\a.e x\in A_{k}$ with $k\in \mathbb{Z}$, applying the size condition (\ref{equ.6}) and the generalized H\"{o}lder's inequality(see (\ref{Hd}) in Lemma \ref{Lem.1}), we have
\begin{equation}\label{equ.10}
\begin{split}
|T_{\beta}(f_{j})(x)|&\le  C2^{j(\beta-n)}\|f_{j}\|_{L^{1}(\R^{n})} \le C2^{j(\beta-n)}\|f_{j}\|_{L^{q_{1}(\cdot)}(\R^{n})}\| \chi_{_{\scriptstyle j}}\|_{L^{q'_{1}(\cdot)}(\R^{n})}.
\end{split}
\end{equation}
Using Proposition \ref{Pr-1}, Proposition \ref{Pr-2}, Proposition \ref{Pr-4}, Lemma \ref{Lem.2}, Lemma \ref{Lem.3} and (\ref{equ.4}), we obtain
\begin{equation}\label{equ.11}
\begin{split}
 2^{j(\beta-n)}\|\chi_{_{\scriptstyle k}}\|_{L^{q_{2}(\cdot)}(\R^{n})} \|\chi_{_{\scriptstyle j}}\|_{L^{q'_{1}(\cdot)}(\R^{n})}
 &\le 2^{j(\beta-n)}\|\chi_{_{B_{k}}}\|_{L^{q_{2}(\cdot)}(\R^{n})}  \|\chi_{_{B_{j}}}\|_{L^{q'_{1}(\cdot)}(\R^{n})} \\
 &\le C\|\chi_{_{B_{k}}}\|_{L^{q_{2}(\cdot)}(\R^{n})} \cdot2^{j\beta}2^{-jn} \|\chi_{_{B_{j}}}\|_{L^{q'_{1}(\cdot)}(\R^{n})}\\
 &\le C \|\chi_{_{B_{k}}}\|_{L^{q_{2}(\cdot)}(\R^{n})} \cdot2^{j\beta}\|\chi_{_{B_{j}}}\|^{-1}_{L^{q_{1}(\cdot)}(\R^{n})} \\
 &\le C \|\chi_{_{B_{k}}}\|_{L^{q_{2}(\cdot)}(\R^{n})} \cdot\|\chi_{_{B_{j}}}\|^{-1}_{L^{q_{2}(\cdot)}(\R^{n})} \\
 &= C \frac{\|\chi_{_{B_{k}}}\|_{L^{q_{2}(\cdot)}(\R^{n})}}{\|\chi_{_{B_{j}}}\|_{L^{q_{2}(\cdot)}(\R^{n})}}\le C 2^{(k-j)n\delta_{4}}.
\end{split}
\end{equation}
Therefore, combining (\ref{equ.9}), (\ref{equ.10}) and (\ref{equ.11}), and using $\alpha> \lambda-n\delta_{4}$, it follows that
\begin{equation*}
\begin{split}
E_{3}&=\sup_{k_{0}\in \mathbb{Z}}2^{-k_{0}\lambda p_{_{1}}}  \bigg(\sum_{k=-\infty}^{k_{0}}2^{k\alpha p_{_{1}}}
  \Big(\sum_{j=k+2}^{\infty}\|T_{\beta}(f_{j})\cdot\chi_{_{\scriptstyle{k}}}  \|_{L^{q_{_{2}}(\cdot)}(\R^{n})} \Big)^{p_{_{1}}}\bigg)\\
&\le C \sup_{k_{0}\in \mathbb{Z}}2^{-k_{0}\lambda p_{_{1}}}  \bigg(\sum_{k=-\infty}^{k_{0}}2^{k\alpha p_{_{1}}}\Big( \sum_{j=k+2}^{\infty} 2^{j(\beta-n)}\|f_{j}\|_{L^{q_{1}(\cdot)}(\R^{n})}\| \chi_{_{\scriptstyle j}}\|_{L^{q'_{1}(\cdot)}(\R^{n})}\|\chi_{_{\scriptstyle{k}}}  \|_{L^{q_{_{2}}(\cdot)}(\R^{n})} \Big)^{p_{_{1}}}\bigg)  \\
&\le C \sup_{k_{0}\in \mathbb{Z}}2^{-k_{0}\lambda p_{_{1}}}  \bigg(\sum_{k=-\infty}^{k_{0}}2^{k\alpha p_{_{1}}}\Big( \sum_{j=k+2}^{\infty} 2^{(k-j)n\delta_{4}}\|f_{j}\|_{L^{q_{1}(\cdot)}(\R^{n})}\Big)^{p_{_{1}}}\bigg)  \\
&\le C\|f\|^{p_{1}}_{M\dot{K}_{p_{1},q_{1}(\cdot)}^{\alpha,\lambda}(\R^{n})}\sup_{k_{0}\in \mathbb{Z}}2^{-k_{0}\lambda p_{_{1}}}  \bigg(\sum_{k=-\infty}^{k_{0}}2^{k\alpha p_{_{1}}}\Big( \sum_{j=k+2}^{\infty} 2^{(k-j)n\delta_{4}}2^{j(\lambda-\alpha)}\Big)^{p_{_{1}}}\bigg)  \\
&\le C\|f\|^{p_{1}}_{M\dot{K}_{p_{1},q_{1}(\cdot)}^{\alpha,\lambda}(\R^{n})}\sup_{k_{0}\in \mathbb{Z}}2^{-k_{0}\lambda p_{_{1}}}  \bigg(\sum_{k=-\infty}^{k_{0}}2^{k\lambda p_{_{1}}}\Big( \sum_{j=k+2}^{\infty} 2^{(k-j)(\alpha-\lambda+n\delta_{4})}\Big)^{p_{_{1}}}\bigg)  \\
&\le C\|f\|^{p_{1}}_{M\dot{K}_{p_{1},q_{1}(\cdot)}^{\alpha,\lambda}(\R^{n})}\sup_{k_{0}\in \mathbb{Z}}2^{-k_{0}\lambda p_{_{1}}}  \Big(\sum_{k=-\infty}^{k_{0}}2^{k\lambda p_{_{1}}} \Big)   \le C\|f\|^{p_{1}}_{M\dot{K}_{p_{1},q_{1}(\cdot)}^{\alpha,\lambda}(\R^{n})}.
\end{split}
\end{equation*}

Combining the estimates for  $E_{1}$,  $E_{2}$ and  $E_{3}$, we conclude that
$$\|T_{\beta}(f)\|_{M\dot{K}_{p_{_{2}}, q_{_{2}}(\cdot)}^{\alpha,\lambda}(\R^{n})} \le C\|f\|_{M\dot{K}_{p_{_{1}}, q_{_{1}}(\cdot)}^{\alpha,\lambda}(\R^{n})}$$
and then completes the proof of  Theorem \ref{thm.2}.
\end{proof}

\begin{remark}
It is easy to see that when $\lambda=0$, the above results are also true on the Herz space with variable exponent, and containing some main results for \cite{I2}.
\end{remark}

\begin{acknowledgments}
 The author cordially  thank the referees for their valuable suggestions and useful comments which have lead to the improvement of this paper. This work was supported   by  the  Project (No. SY201313) of Mudanjiang Normal University and the Research Project (No.12531720) for Department of Education  of Heilongjiang Province.
\end{acknowledgments}

\end{document}